\documentclass{ifacconf}

\usepackage{graphicx}      % include this line if your document contains figures
\usepackage[table]{xcolor}

\makeatletter
\let\old@ssect\@ssect % Store how ifacconf defines \@ssect
\makeatother
\usepackage{natbib}
\usepackage[hidelinks]{hyperref}
\makeatletter
\def\@ssect#1#2#3#4#5#6{%
  \NR@gettitle{#6}% Insert key \nameref title grab
  \old@ssect{#1}{#2}{#3}{#4}{#5}{#6}% Restore ifacconf's \@ssect
}
\makeatother

\usepackage{mathtools}
\usepackage{amsmath}
\usepackage{amssymb}
\usepackage{enumitem}
\usepackage{booktabs}  
\usepackage{arydshln} %draw dashed lines in arrays
\usepackage{multirow}
\definecolor{Gray}{gray}{0.9}
\usepackage{float}
\usepackage{anyfontsize}

\newcommand{\Rank}[1]{\mathrm{rank}\left(#1\right)}
\newcommand{\Dim}[1]{\mathrm{dim}(#1)}

\newcommand{\D}{\mathrm{d}}

\newcommand{\Qm}{\mathcal{Q}}

\begin{document}
\begin{frontmatter}

\title{ On the Exact Linearization of Minimally Underactuated Configuration Flat Lagrangian Systems in \\ Generalized State Representation } 

\thanks[footnoteinfo]{This research was funded in whole, or in part, by the Austrian Science Fund (FWF) P36473. For the purpose of open access, the author has applied a CC BY public copyright licence to any Author Accepted Manuscript version arising from this submission.}

\author[First]{Georg Hartl}
\author[First]{Conrad Gstöttner}
\author[First]{Bernd Kolar}
\author[First]{Markus Schöberl}

\address[First]{Institute of Automatic Control and Control Systems Technology, Johannes Kepler University Linz, Altenberger Strasse 69, 4040 Linz, Austria (e-mail: \{georg.hartl, conrad.gstoettner, bernd.kolar, markus.schoeberl\}@jku.at).}

\begin{abstract} % Abstract of 50--100 words

In this paper, we examine the exact linearization of configuration flat Lagrangian control systems in generalized state representation with $p$ degrees of freedom and $p-1$ control inputs by quasi-static feedback of its generalized state. We formally introduce generalized Lagrangian control systems, which are obtained when configuration variables are considered as inputs instead of forces. This work presents all possible lengths of integrator chains achieved by an exact linearization with a quasi-static feedback law of the generalized state that allows for rest-to-rest transitions. We show that such feedback laws can be systematically derived without using Brunovský states.\vspace{-10pt}
\end{abstract}
\begin{keyword}
Flatness, Lagrangian Control Systems, Differential Geometry 
\end{keyword}
\end{frontmatter}
%===============================================================================

\section{ Introduction }
\vspace{-8pt}
The concept of differential flatness, as introduced, e.g., in \cite{FliessFlatnessDefectNonlinear1995}, has significantly advanced nonlinear control theory. Roughly speaking, a nonlinear control system in classical state representation $\dot{x} = f(x,u)$ is flat if the $n$-dimensional classical state $x$ and the $m$-dimensional control input $u$ can be expressed in terms of an $m$-tuple $y$ of differentially independent functions $y_j=\varphi_j(x,u,u^{(1)}, \ldots, u^{(\nu)})$, $j=1,\ldots,m$, and their time derivatives, where $u^{(\nu)}$ denotes the $\nu$-th time derivative of the input. However, computing a so-called flat output $y$ continues to be a complex challenge lacking a comprehensive solution, as emphasized in recent studies such as \cite{GstottnerFlatSystemPossessing2023}, \cite{GstottnerNecessarySufficientConditions2023}, \cite{NicolauNormalFormsMulti2019}, \cite{NicolauFlatnessMultiInputControlAffine2017}, and \cite{SchoberlImplicitTriangularDecomposition2014}. Exploiting the property of flatness enables elegant and systematic trajectory planning and tracking solutions. The concept of flatness-based tracking control design typically involves exact linearization, where a nonlinear control system is exactly linearized by endogenous feedback, resulting in a closed-loop system consisting of $m$ integrator chains \mbox{$y_j^{(\kappa_j)} = w_j$}, \mbox{$j=1,\ldots,m$}, with respective lengths $\kappa_j$ between a new input $w$ and the flat output $y$. As shown, e.g., in \cite{FliessLieBacklundApproachEquivalence1999}, the exact linearization of flat systems is systematically achievable by applying dynamic feedback laws. However, this approach has the disadvantage of increasing the order of the closed-loop system dynamics, consequently resulting in an increased order of the tracking error dynamics. As an alternative, \cite{DelaleauControlFlatSystems1998} show that every flat system is linearizable by a quasi-static feedback law \mbox{$u=\alpha(\Tilde{x}_B, w, w^{(1)}, \ldots)$} of a so-called generalized Brunovský state $\Tilde{x}_B$ that consists of suitably chosen consecutive time derivatives of the flat output $y$ meeting $\Dim{\Tilde{x}_B}=\Dim{x}$. Consequently, applying such a feedback law preserves the order of the closed-loop system dynamics. However, practically measuring or determining $\Tilde{x}_B$ from available measurement data is often rather complex. Hence, the exact linearization of nonlinear control systems by a quasi-static feedback law $u=\alpha(x,w,w^{(1)}, \ldots )$ of the classical state $x$ is of significant interest. For systems in classical state representation with a flat output of the form $ y = \varphi(x, u) $, known as $(x, u)$-flat systems, \cite{DelaleauFiltrationsFeedbackSynthesis1998} demonstrate in a differential algebraic setting that a linearizing feedback law $u = \alpha(x, w, w^{(1)}, \ldots)$ is guaranteed to exist. In a differential geometric framework, \cite{GstottnerTrackingControlFlat2024} derive similar results and present an algorithm that constructs the lengths of the integrator chains of the closed-loop system using successive coordinate transformations, avoiding Brunovský states. \cite{KolarExactLinearisationControl2022} extend these results to differentially flat discrete-time systems.

Since various mechanical systems are flat, the analysis of flatness of mechanical control systems is an extensively researched area, as highlighted, e.g., in \cite{NicolauFlatnessMechanicalSystems2015} and \cite{KnollConfigurationFlatnessLinear2014}. Yet, within this paper, we restrict ourselves to flat systems characterized by Lagrangian mechanics with a Lagrangian of the form kinetic minus potential energy. 
% The control input of Lagrangian control (LC) systems in classical state representation is given by forces and torques with the corresponding classical state consisting of the configuration variables $q$ and their respective velocities $v$. 
Lagrangian control (LC) systems with $p$ degrees of freedom and $p-1$ control inputs that possess a flat output depending only on configuration variables are called minimally underactuated configuration flat LC systems, see, e.g., \cite{RathinamConfigurationFlatnessLagrangian1998}. A detailed analysis of the exact linearization of this class of control systems by quasi-static feedback of the classical state is presented in \cite{HartlExactLinearizationMinimally2023a}. While classical LC systems are predominantly addressed in both practical applications and academic literature, considering LC systems in generalized state representation reveals distinct advantages. 

Control systems in generalized state representation involve system dynamics that depend on the generalized state $\Tilde{x}$, the control input $u$, and a finite number of time derivatives of $u$. For a more detailed discussion, we refer, e.g., to \cite{FliessGeneralizedControllerCanonical1990a}. As presented, e.g., in \cite{FliessSimplifiedApproachCrane1991} or \cite{RudolphExamplesRemarksQuasiStatic1998}, generalized LC systems arise naturally when positions are considered as control inputs and the corresponding forces are determined by underlying control loops using, e.g., PI, PD, or PID controllers.\footnote{In our work, we assume that the underlying control loops operate sufficiently well so that they do not need further consideration.} The resulting system dynamics involve time derivatives of the newly introduced control inputs. As demonstrated, e.g., in \cite{FliessSimplifiedApproachCrane1991}, generalized state representations enable simplified physical models and prove beneficial for designing more robust control laws.

This paper discusses the exact linearization of minimally underactuated configuration flat generalized LC systems by quasi-static feedback of its generalized states. These generalized states consist of those configuration variables $\Tilde{q}$ that are not selected as control inputs and their respective velocities $\Tilde{v}$. Our approach builds upon essential results outlined in \cite{GstottnerTrackingControlFlat2024} and \cite{HartlExactLinearizationMinimally2023a}, and extends them to generalized LC systems. Based on the Jacobian of the parameterization of the system variables $(\Tilde{q},\Tilde{v})$ by the flat output $y$, we demonstrate how to choose suitable orders of the time derivatives of the flat output intended as new inputs. We present all possible lengths of integrator chains that can be achieved by applying a quasi-static feedback law of the generalized state $(\Tilde{q},\Tilde{v})$ that exhibits no singularities at equilibrium points. Our approach uses neither successive coordinate transformations nor Brunovský states. The structure of this paper is organized as follows: Section \ref{sec_2} establishes the notation. Section \ref{sec_3} examines properties of differentially flat minimally underactuated classical LC systems. Section \ref{sec_4} covers the principles of differential flatness and exact linearization for generalized state systems. Section \ref{sec_5} formally introduces generalized LC systems. Finally, Section \ref{sec_6} is devoted to our main results, accompanied by a practical example.
\vspace{-5pt}
\section{Notation}
\label{sec_2}
\vspace{-5pt}
This work utilizes tensor notation and the Einstein summation convention. Consider an $n$-dimensional smooth manifold $\mathcal{M}$ with local coordinates $x=(x^1, \ldots, x^n)$ and an $m$-tuple of functions $h=(h^1, \ldots, h^m):\mathcal{M} \rightarrow \mathbb{R}^m$ defined on $\mathcal{M}$. By $\partial_x h$, we denote the $m \times n$ Jacobian matrix, and $\partial_{x^i}h^j$ represents the partial derivative of $h^j$ with respect to $x^i$. The notation $\D h$ represents the differentials $(\D h^1, \ldots, \D h^m)$. We denote time derivatives using subscripts in square brackets. For example, $y^j_{[\alpha]}$ represents the $\alpha$-th time derivative of the $j$-th element within the $m$-tuple $y$, and $y_{[\alpha]} = (y^1_{[\alpha]}, \dots , y^m_{[\alpha]})$ indicates the $\alpha$-th time derivative for each component of  $y$. Capitalized multi-indices denote time derivatives of different orders of each element in a tuple. Given two multi-indices, $A = (a^1, \dots, a^m)$ and $B=(b^1, \dots, b^m)$, we can find concise representations of time derivatives of different orders, such as \mbox{$y_{[A]} = (y^1_{[a^1]}, \dots,y^m_{[a^m]})$}, \mbox{$y_{[0, A]} = ( y^1_{[0,a^1]}, \dots,y^m_{[0,a^m]} )$}, and \mbox{$y_{[A, B]} = ( y^1_{[a^1,b^1]}, \dots,y^m_{[a^m,b^m]} )$}. For \( a^j \leq b^j \), the sequence \( y^j_{[a^j, b^j]} \) denotes successive time derivatives and is defined as \( y^j_{[a^j, b^j]} = ( y^j_{[a^j]}, \dots, y^j_{[b^j]}) \). Conversely,
\begin{equation}\label{conv_1}
    a^j > b^j \implies y^j_{[a^j, b^j]} = \emptyset.
\end{equation}
We express the component-wise addition and subtraction of multi-indices as $A \pm B = ( a^1\pm b^1, \dots, a^m \pm b^m )$. The summation over the indices is denoted as $\#A = \sum_{j=1}^m a^j$.

\vspace{-6pt}
\section{Classical Lagrangian Control Systems}
\label{sec_3}
\vspace{-6pt}

In this section, we briefly review classical LC systems and refine essential findings from \cite{HartlExactLinearizationMinimally2023a}. Consider the $p$-dimensional configuration space $\Qm$ with $(q^1, \ldots, q^p)$. As in \cite{RathinamConfigurationFlatnessLagrangian1998}, we focus on mechanical systems whose dynamics can be described by employing the Lagrangian formalism with a Lagrangian 
\begin{equation*}
    L(q,\dot{q})=\frac{1}{2}g(\dot{q}, \dot{q})-V(q)
\end{equation*} 
of the form kinetic minus potential energy, where $g$ denotes a Riemannian metric on $\Qm$. Further, we assume that the implicit equations of motion are given by
\begin{equation}
    \label{equ_of_motion}
    \frac{\D}{\D t} \left( \frac{\partial L(q,\dot{q})}{\partial \dot{q}^i} \right) - \frac{\partial L(q,\dot{q})}{\partial q^i} = G_{ij}(q)u^j,
\end{equation}
with $i=1,\ldots,p$ and $j=1,\ldots,m$. Solving the equation system (\ref{equ_of_motion}) for  $\ddot{q}$ yields the classical state representation
\begin{subequations}
    \label{lag_sys_class}
    \begin{align}
        \dot{q}^i & = v^i,  \label{lag_sys_class_q} \\
        \dot{v}^i & = f^i(q,v,u), \quad i=1,\ldots,p, \label{lag_sys_class_v} 
    \end{align}
\end{subequations}
with the $2p$-dimensional classical state $(q,v)$ and the \mbox{$m$}-dimensional input $u$. To proceed, we concisely review the key findings of \cite{HartlExactLinearizationMinimally2023a} relevant to this paper. For a LC system (\ref{lag_sys_class}) with a flat output $y=(q,v,u,u_{[1]}, \ldots)$, the parameterization of the system variables $(q,v,u)$ by the flat output $y$ is given by 
\begin{subequations}\label{para_map_lag}
    \begin{align}
        (q,v) &= (F_q(y_{[0,R-2]}), F_v(y_{[0,R-1]})), \label{para_map_lag_q_v}\\
        u & = F_u(y_{[0,R]}), \label{para_map_lag_u} \vspace{-1mm}
    \end{align}
\end{subequations} 
with the multi-index $R=(r^1,\ldots,r^m)$ denoting the highest orders of time derivatives that explicitly occur in (\ref{para_map_lag}). 
For LC systems with $m=p-1$ control inputs and a configuration flat output
\begin{equation}\label{class_config_flat_output}
    y^j=\varphi^j(q), \quad j=1,\ldots,p-1,
\end{equation}
$R=(r^1,\ldots,r^{p-1})$ satisfies $2\leq R \leq 4$ with
\begin{equation}
    \label{R_1}
    r^j = 4 \text{ for at least one } j \in \{ 1,\ldots,p-1 \}.
\end{equation} 
Further, there always exists a regular state transformation \vspace{-4mm}
\begin{subequations}\label{state_trf_1}
    \begin{align}
        \bar{q}^j & = g^j(q) = \varphi^j(q), \quad j=1,\ldots,p-1, \label{state_trf_1_a} \\
        \bar{q}^p & = g^p(q), \label{state_trf_1_b} \\
        \bar{v}^i & = \partial_{q^l} g^i(q)v^l, \hspace{9.4mm} i,l=1,\ldots,p. \label{state_trf_1_c}
    \end{align}
\end{subequations}
Consequently, the parameterization of the transformed coordinates $(\bar{q}, \bar{v})$ by the flat output is of the form
\begin{equation}
    \label{flat_para_sec_2}
    \arraycolsep=3pt
    \begin{array}{rclcrcl}
        \bar{q}^1 &=& y^1, & & \bar{v}^1 &=& y^1_{[1]}, \\ [-0.1cm]
        &\vdots&  & &  &\vdots&  \\ [-0.1cm]
        \bar{q}^{p-1} &=& y^{p-1}, & & \bar{v}^{p-1} &=& y^{p-1}_{[1]}, \\ [0.2cm]
        \bar{q}^{p} &=& F^p_{\bar{q}} \left( y_{[0, R-2 ]} \right), & & \bar{v}^{p} &=& F^p_{\bar{v}} \left( y_{[0, R-1 ]} \right). \\ [0.2cm]
    \end{array}
        \vspace{-3pt}
\end{equation}
Further, Appendix A of \cite{HartlExactLinearizationMinimally2023a} shows that the parameterization of $\bar{q}^p$ has the form $F_{\bar{q}}^p(y, y^l_{[1]}y^k_{[1]}, y_{[2]})$. In other words, if $F_{\bar{q}}$ explicitly depends on first-order time derivatives of $y$, they can only occur quadratically. 
 
In the derivation of feedback laws that enable rest-to-rest transitions, equilibrium points are particularly important. For flat LC systems, equilibrium points $(q_s,0,u_s)$ are parameterized by constant values $y_s$ of the flat output $y$, with $y_{s,[1]}=y_{s,[2]}=\ldots=0$. In this work, we assume that the flat parameterization (\ref{para_map_lag}) is well defined around any considered equilibrium point $y_s$.\footnote{Note that \textit{around} refers to a neighborhood of $y_s$ on the smooth manifold $\mathcal{Y}_{[0,R]}$ with local coordinates $(y,\ldots,y_{[R]})$. } The following lemma refines core results of \cite{HartlExactLinearizationMinimally2023a} by demonstrating that the Jacobian matrix of $(F_q, F_v)$ evaluated at equilibrium points exhibits a specific structure.
% Given that $F_q$ depends quadratically on first-order time derivatives of $y$, the following lemma can be easily derived. 
\begin{lem}\label{lem_2}
    Consider a minimally underactuated LC system (\ref{lag_sys_class}) with a configuration flat output (\ref{class_config_flat_output}) as well as the flat parameterization (\ref{para_map_lag}). The Jacobian of $(F_{q}, F_{v})$ evaluated at equilibrium points is given by
    \begin{equation}
        \label{jacobian_1_ys_classic}
        \begin{aligned}
        \left.
        \partial_{y_{[0,3]}} 
        \begin{pmatrix}
            F_{q} \\ F_{v}
        \end{pmatrix}\right|_{y_s} \hspace{-5pt} = 
        \left( \begin{array}{cccc}
            \partial_{y}F_{q}|_{y_{s}} & 0 & \partial_{y_{[2]}}F_{q}|_{y_{s}} & 0 \\
            0 & \partial_{y}F_{q}|_{y_{s}} & 0 & \partial_{y_{[2]}}F_{q}|_{y_{s}} \\
        \end{array} \right),
        \end{aligned} %\vspace{12pt}
    \end{equation} %we assume that (\ref{para_map_lag}) is well-defined locally around any considered equilibrium point \((q_s, 0, u_s)\) parameterized by \(y_s\),
    with\footnote{ We assume that the rank of \(\partial_{y_{[2]}}F_q\) remains constant around any considered $y_s$. } \vspace{-3mm}
    % In this work, we assume that the flat parameterization (\ref{para_map_lag}) is well-defined and that the rank of \(\partial_{y_{[2]}}F_q\) remains constant, both locally around any considered equilibrium point \((q_s, 0, u_s)\) parameterized by \(y_s\).} \vspace{-3mm}
    \begin{subequations}\label{lem_1_rank}
        \begin{align}
            \Rank{\partial_y F_{q}|_{y_s}} &= p-1, \label{lem_1_rank_a} \\
            \Rank{\partial_{y_{[2]}}F_{q}|_{y_s}} &\leq 1. \label{lem_1_rank_b}
        \end{align}
    \end{subequations}
\end{lem}
\begin{pf}
    Given $2\leq R \leq 4$, calculating the Jacobian of (\ref{para_map_lag_q_v}) only needs the consideration of time derivatives of $y$ up to the third order. First, we analyze the upper $p$ rows of (\ref{jacobian_1_ys_classic}). Recall that $F_{\bar{q}}^p$ in (\ref{flat_para_sec_2}) depends quadratically on first-order time derivatives of $y$. From (\ref{state_trf_1_a})-(\ref{state_trf_1_b}), it follows that $F_q$ also depends quadratically on first-order time derivatives of $y$, i.e., $q = F_q(y, y^l_{[1]}y^k_{[1]}, y_{[2]})$. Hence, $\partial_{y_{[1]}} F_q \big|_{y_s}$ and $\partial_{y_{[3]}} F_q$ are zero.
    Next, we analyze the lower $p$ rows of (\ref{jacobian_1_ys_classic}). Since the relationship between $F_v$ and $F_q$ is given by
    \begin{equation*}
        F^i_v = y^j_{[1]}\partial_{y^j}F^i_q + y^j_{[2]}\partial_{y^j_{[1]}}F^i_q + y^j_{[3]}\partial_{y^j_{[2]}}F^i_q,
    \end{equation*}
    with $i = 1,\ldots,p,$ and $j=1,\ldots,p-1$, it follows that
    \begin{equation*}\label{plan_manip_class}
    \arraycolsep=3pt
    \begin{array}{rclcrcl}
        \partial_{y^j}F^i_v\big|_{y_{s}} & = & 0, & & \partial_{y^j_{[2]}}F^i_v\big|_{y_{s}} & = & \partial_{y^j_{[1]}}F^i_q\big|_{y_{s}} = 0, \\ [0.2cm]
        \partial_{y^j_{[1]}}F^i_v\big|_{y_{s}} & = & \partial_{y^j}F^i_q\big|_{y_{s}}, & & \partial_{y^j_{[3]}}F^i_v\big|_{y_{s}} & = & \partial_{y^j_{[2]}}F^i_q\big|_{y_{s}}. \\ [0.2cm]
    \end{array}
\end{equation*}
Next, let us prove (\ref{lem_1_rank}). By examining the left column of (\ref{flat_para_sec_2}) it immediately follows that $\text{rank}(\partial_y F_{\bar{q}}|_{y_s}) = p-1$. Because of (\ref{R_1}), $F^p_{\bar{q}}$ depends on the second-order time derivative of at least one component of $y$. Therefore, we find that $\text{rank}(\partial_{y_{[2]}} F_{\bar{q}}) = 1$, which implies that $\text{rank}(\partial_{y_{[2]}} F_{\bar{q}}|_{y_s}) \leq 1$. Finally, the diffeomorphism (\ref{state_trf_1_a})-(\ref{state_trf_1_b}) implies that the same rank conditions also apply to the Jacobian matrices $\partial_yF_q|_{y_s}$ and $\partial_{y_{[2]}}F_q|_{y_s}$ of the parameterization $F_q$ of the original coordinates $q$. \hfill $\Box$
\end{pf}
\vspace{-5pt}
\section{ Flatness of Generalized Systems }\label{sec_4}
\vspace{-5pt}
In this section, we consider nonlinear control systems in generalized state representation
\begin{equation}\label{gen_sys_1}
    \dot{\Tilde{x}}=\Tilde{f}( \Tilde{x}, u_{[0,B]} ), \\
\end{equation}
with the $n$-dimensional state $\Tilde{x}$, the $m$-dimensional input $u$ and the multi-index $B=(b^1, \ldots, b^m)$, $B \geq 0$, denoting the highest orders of time derivatives of the input components occurring in the right-hand side of (\ref{gen_sys_1}). First, we present a definition of differentially flat generalized state systems (\ref{gen_sys_1}) by utilizing a finite-dimensional differential geometric framework as, e.g., in \cite{KolarPropertiesFlatSystems2016}. We use a state-input manifold $\Tilde{\mathcal{X}}\times\mathcal{U}_{[0,l_u]}$ with local coordinates $(\Tilde{x}, u, u_{[1]}, \ldots, u_{[l_u]})$, including time derivatives of the input $u$ up to some finite but large enough order $l_u$.
\begin{defn}\label{def_2}
    A nonlinear system (\ref{gen_sys_1}) is called flat if there exists an $m$-tuple of smooth functions $y^j=\varphi^j(\Tilde{x}, u, u_{[1]}, \ldots )$, $j=1,\ldots,m$, defined on $\Tilde{\mathcal{X}}\times\mathcal{U}_{[0,l_u]}$ and smooth functions $F^i_{\Tilde{x}}$ and $F^j_{u}$ such that locally
    \begin{subequations}\label{para_map_def}
        \begin{align}
            \Tilde{x}^i & = F^i_{\Tilde{x}}(y_{[0,R-1]}), \hspace{4mm} i=1,\ldots,n, \label{para_map_def_x} \\ 
            u^j & = F^j_{u}(y_{[0,S]}), \hspace{7.4mm} j=1,\ldots,m, \label{para_map_def_u}
        \end{align}
    \end{subequations}
    with multi-indices $R=(r^1, \ldots, r^m)$, $S=(s^1, \ldots, s^m)$, and $R, S \geq 0$. The $m$-tuple $y$ is called a flat output of (\ref{gen_sys_1}), and (\ref{para_map_def}) is called the parameterizing map.
\end{defn}
The objective of exact linearization is to establish a linear input-output behavior $y^j_{[\kappa^j]} = w^j, j=1,\ldots,m,$ between a newly introduced input $w$ and the flat output $y$ by applying a suitable feedback law. However, a major challenge lies in choosing the appropriate lengths $\kappa$ of the integrator chains. In the context of tracking-control design, it is beneficial to reduce the orders of the tracking error dynamics by selecting the orders $\kappa^j$ as low as possible, which can be achieved by constructing a linearizing quasi-static feedback law of the generalized state $\Tilde{x}$. In the following, we present conditions that ensure that selected time derivatives of a flat output can be introduced as new closed-loop input $w$ by applying a quasi-static feedback law of $\Tilde{x}$. With the focus on quasi-static feedback of $\Tilde{x}$, we extend Theorem 3.1 from \cite{GstottnerTrackingControlFlat2024} to generalized state systems (\ref{gen_sys_1}).
\begin{thm}\label{thm_3}
    Consider system (\ref{gen_sys_1}) with a flat output $y=\varphi(\Tilde{x},u,u_{[1]}, \ldots)$. An $m$-tuple of time derivatives of this flat output given by $y_{[\kappa]}=\varphi_{[\kappa]}(\Tilde{x},u,u_{[1]}, \ldots)$, with $\kappa \leq R$ and $\# \kappa = n$, can be introduced as a new input $w$ by a quasi-static feedback law $u = \alpha(\Tilde{x}, w, w_{[1]},\ldots)$ of the generalized state $\Tilde{x}$ if and only if the differentials $\D \Tilde{x}, \D \varphi_{[\kappa]}, \D \varphi_{[\kappa+1]}, \ldots$ are linearly independent.\footnote{The proof is omitted here for brevity but can be conducted analogously to the proof of Theorem 3.1 in \cite{GstottnerTrackingControlFlat2024}.}
\end{thm}
\vspace{-7pt}
\section{ Generalized LC systems }\label{sec_5}
\vspace{-7pt}
The starting point of this section is the classical state representation (\ref{lag_sys_class}) of an underactuated LC system (i.e., $m<p$). We aim to introduce $1 \leq k \leq m$ reasonably chosen configuration variables as new control inputs. In general, we consider replacing $k$ classical control input components $(u^1, \ldots, u^k)$ by $k$ configuration variables $(q^{1}, \ldots, q^{k})$. In other words, we assume that the forces and torques $(u^1, \ldots, u^k)$, which act as control inputs of the classical system, are well determined by low-level controllers and are not further involved in the system dynamics. Given a system in classical representation (\ref{lag_sys_class}), the generalized system dynamics are obtained by solving (\ref{lag_sys_class_v}) for
\begin{subequations}\label{sol_gen_sys}
    \begin{align}
        u^i & = h^i\left( q, v, \dot{v}^{1}, \ldots, \dot{v}^{k}, u^{k+1},\ldots,u^{m} \right), \label{sol_gen_sys_a} \\
    \dot{v}^j & = h^j\left(q, v, \dot{v}^{1}, \ldots, \dot{v}^{k}, u^{k+1},\ldots,u^{m}\right), \label{sol_gen_sys_b}
    \end{align}
\end{subequations}
with $i=1,\ldots,k,$ and $j=k+1,\ldots,p$. To distinguish generalized from classical LC systems, we mark variables and indices associated with generalized LC systems with the tilde symbol and substitute 
\begin{subequations}\label{gen_subs}
    \begin{align}
    (q^i, v^i, \dot{v}^i) &= (\Tilde{u}^i_{[0,2]}), \hspace{18.5mm} i=1,\ldots,k, \label{gen_subs_a} \\
        u^l & = \Tilde{u}^l, \hspace{25.5mm} l=k+1,\ldots,m, \label{gen_subs_b} \\
        (q^j, v^j, \dot{v}^j) &= (\Tilde{q}^{j-k}, \Tilde{v}^{j-k}, \dot{\Tilde{v}}^{j-k}), \hspace{1mm} j=k+1,\ldots,p \label{gen_subs_c}
    \end{align}
\end{subequations} 
into (\ref{sol_gen_sys_b}). Hence, a generalized LC system is of the form
\begin{equation}
    \label{lag_gen_sys_1}
    \begin{aligned}
        \dot{\Tilde{q}}^i & = \Tilde{v}^i, \\
        \dot{\Tilde{v}}^i & = \Tilde{f}^i ( \Tilde{q}, \Tilde{v}, \Tilde{u}_{[0,\Tilde{B}]} ), \quad i=1,\ldots,p-k,
    \end{aligned}
\end{equation}
with the $2(p-k)$-dimensional generalized state $(\Tilde{q}, \Tilde{v})$, the input $\Tilde{u}$, and the multi-index $\Tilde{B} = (\Tilde{b}^1, \ldots, \Tilde{b}^k, 0, \ldots, 0)$, with $0 \leq \Tilde{B} \leq 2$. Throughout this work, we assume that all considered generalized state representations (\ref{lag_gen_sys_1}) are well-defined locally around equilibrium points. Let us now consider a classical LC system (\ref{lag_sys_class}) with a flat output $y=(q, v, u, u_{[1]}, \ldots)$ permitting a generalized representation (\ref{lag_gen_sys_1}). The resulting generalized LC system is also flat, and the parameterization of the system variables $(\Tilde{q}, \Tilde{v}, \Tilde{u})$ reads 
\begin{subequations}\vspace{-4mm}
    \label{para_map_gen}
    \begin{align}
        \label{para_map_gen_q} (\Tilde{q},\Tilde{v}) & = (F_{\Tilde{q}}(y_{[0,\Tilde{R}-2]}), F_{\Tilde{v}}(y_{[0,\Tilde{R}-1]})),  \\
        \label{para_map_gen_u} \Tilde{u} & = F_{\Tilde{u}}(y_{[0,\Tilde{S}]}), 
    \end{align}
\end{subequations}
with the multi-indices $\Tilde{R}=(\Tilde{r}^1, \ldots, \Tilde{r}^m)$ and $\Tilde{S}=(\Tilde{s}^1, \ldots, \Tilde{s}^m)$. According to (\ref{gen_subs_a})-(\ref{gen_subs_c}), the generalized state variables are $(\Tilde{q}, \Tilde{v})=(q^{k+1}, \ldots, q^{p}, v^{k+1}, \ldots, v^{p})$ and the control input reads $\Tilde{u}=(q^1,\ldots,q^k,u^{k+1},\ldots,u^m)$.\footnote{ Our approach is built upon the assumption that measurements or estimates of the classical state are available. Given that the generalized state $(\Tilde{q}, \Tilde{v})$ comprises only components of $(q,v)$, it follows that measurements or estimates of $(\Tilde{q}, \Tilde{v})$ are likewise known. } Hence, the flat parameterization (\ref{para_map_gen}) is simply obtained from (\ref{para_map_lag}) by omitting the parameterizations of $(u^1, \ldots, u^k)$ and $(v^{1}, \ldots, v^{k})$. For the special case that the original classical LC system possesses a configuration flat output (\ref{class_config_flat_output}), it follows that (\ref{lag_gen_sys_1}) also possesses a flat output
\begin{equation}\label{flat_output_gen_sys}
    y^j = \varphi^j(\Tilde{q}, \Tilde{u}^1, \ldots, \Tilde{u}^k), \quad j=1,\ldots,m.
\end{equation}
Hence, we still speak of a configuration flat output and, thus, of a configuration flat generalized LC system. Although we are examining generalized state representations, we only consider regular state transformations of the form
\begin{equation}
    \begin{aligned}
        \bar{\Tilde{q}}^i & = g^i( \Tilde{q} ), \\
        \bar{\Tilde{v}}^i & = \partial_{\Tilde{q}^j} g^i(\Tilde{q})\Tilde{v}^j, \quad i,j = 1, \ldots, p-k.
    \end{aligned}
\end{equation}
\vspace{-20pt}
\begin{figure}[H]
    \centering
    \def\svgwidth{0.35\textwidth}
    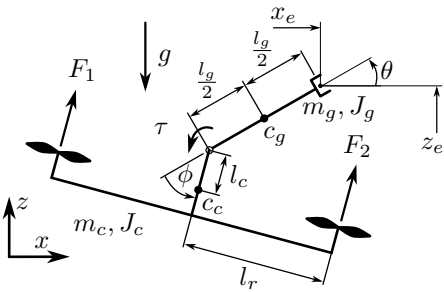
    \caption{Planar aerial manipulator}
    \label{fig1} 
\end{figure}
\vspace{-3mm}
\subsection{ An example illustrating generalized LC systems }\label{sec_lag_con}
\vspace{-8pt}
Consider the planar aerial manipulator with $p=4$ degrees of freedom and $m=3$ inputs as shown in Figure \ref{fig1}. The classical representation (\ref{lag_sys_class}) of the aerial manipulator has the classical state $(q,v)=(x_e, z_e, \theta, \phi, v_{x_e}, v_{z_e}, \omega_\theta, \omega_{\phi})$ and the control input $u=(F_1, F_2, \tau)$, where $F_1$ and $F_2$ are the forces generated by the rotors, and $\tau$ is the torque acting on the joint. As shown, e.g., in \cite{WeldeRoleSymmetryConstructing2023}, a configuration flat output of the aerial manipulator reads\vspace{0.5mm}
\begin{equation}\label{plan_man_flat_output}
    (y^1, y^2, y^3) = (x_e - r_e\cos(\phi+\theta), z_e - r_e\sin(\phi+\theta), \theta),\vspace{0.5mm}
\end{equation}
with $r_e = l_c\frac{m_c}{m_g+m_c}$. Assuming that it is advantageous to control the joint angle $\phi$ through underlying control loops, we replace the input component $\tau$ with the configuration variable $\phi$ (i.e., $k=1$). By solving the corresponding equation system (\ref{lag_sys_class_v}) for $(\tau, \dot{v}_{x_e}, \dot{v}_{z_e}, \dot{\omega}_\theta)$ and substituting $(\phi, \omega_\phi, \dot{\omega}_\phi)=(\phi_{[0,2]})$ into the solution of the form (\ref{sol_gen_sys_b}), we obtain a generalized LC system of the form\vspace{-0.5mm}
\begin{equation}\label{plan_manip_gen}
    \arraycolsep=3pt
    \begin{array}{rclcrcl}
        \dot{x}_e & = & v_{x_e}, & & \dot{v}_{x_e} & = & \Tilde{f}_{v_{x_e}}( \theta, \omega_\theta, \phi_{[0,2]}, F_1, F_2 ), \\ [0.2cm]
        \dot{z}_e & = & v_{z_e}, & & \dot{v}_{z_e} & = & \Tilde{f}_{v_{z_e}}( \theta, \omega_\theta, \phi_{[0,2]}, F_1, F_2 ), \\ [0.2cm]
        \dot{\theta} & = & \omega_\theta, & & \dot{\omega}_\theta & = & \Tilde{f}_{\omega_\theta}( \omega_\theta, \phi_{[0,2]}, F_1, F_2 ), \\
    \end{array}\vspace{-0.5mm}
\end{equation}
with the $2(p-1)=6$-dimensional generalized state $(\Tilde{q},\Tilde{v})=(x_e, z_e, \theta, v_{x_e}, v_{z_e}, \omega_\theta)$ and the control input $\Tilde{u}=(\phi, F_1, F_2)$. The resulting system is also configuration flat with the flat output (\ref{plan_man_flat_output}), the parameterizing map
\begin{subequations}\label{plan_manip_flat_para}
    \begin{align}
        \Tilde{q} & = F_{\Tilde{q}}( y^1, y^1_{[2]}, y^2, y^2_{[2]}, y^3_{[0,2]} ), \label{plan_manip_flat_para_a} \\
        \Tilde{v} & = F_{\Tilde{v}}( y^1_{[1,3]}, y^2_{[1,3]}, y^3_{[0,3]} ), \label{plan_manip_flat_para_b} \\
        \Tilde{u} & = F_{\Tilde{u}}( y^1_{[2,4]}, y^2_{[2,4]}, y^3_{[0,4]} ), \label{plan_manip_flat_para_c}
    \end{align}
\end{subequations}
and the multi-indices $\Tilde{R}=(4,4,4)$ and $\Tilde{S}=(4,4,4)$.
\vspace{-4pt}
\section{ Exact Linearization }\label{sec_6}
\vspace{-4pt}
In this section, we examine the exact linearization of generalized LC systems (\ref{lag_gen_sys_1}) with $m=p-1$ inputs and a configuration flat output (\ref{flat_output_gen_sys}). By applying a quasi-static feedback law 
\begin{equation}\label{quasistat_fb}
    \Tilde{u} = \Tilde{\alpha}(\Tilde{q}, \Tilde{v}, \Tilde{w}, \Tilde{w}_{[1]}, \ldots)
\end{equation}
of the state $(\Tilde{q}, \Tilde{v})$, the closed-loop system
\begin{equation}\label{closed_loop}
    \begin{aligned}
        \dot{\Tilde{q}} & = \Tilde{v}, \\
        \dot{\Tilde{v}} & = \Tilde{f}(\Tilde{q}, \Tilde{v},\Tilde{\alpha}_{[0,2]}(\Tilde{q}, \Tilde{v}, \Tilde{w}, \Tilde{w}_{[1]}, \ldots) )
    \end{aligned}
\end{equation}
should exhibit a linear input-output behavior
\begin{equation}\label{input_output_gen}
    y^j_{[\Tilde{\kappa}^j]}=\Tilde{w}^j, \quad j=1,\ldots,p-1,
\end{equation}
between the new input $\Tilde{w}$ and the flat output $y$. In the following, we outline the process for selecting appropriate orders $\Tilde{\kappa}$ of the time derivatives of $y$. But first, we reformulate Theorem \ref{thm_3} for generalized LC systems (\ref{lag_gen_sys_1}) and simplify the condition that the differentials $\D \Tilde{q}, \D \Tilde{v}, \D \varphi_{[\Tilde{\kappa}]}, \D \varphi_{[\Tilde{\kappa}+1]}, \ldots$ must be linearly independent to a rank condition.
\begin{prop}\label{prop_4}
    Consider a generalized LC system (\ref{lag_gen_sys_1}) with a flat output $y=\varphi(\Tilde{q},\Tilde{v},\Tilde{u},\Tilde{u}_{[1]}, \ldots)$ and the parameterizing map (\ref{para_map_gen}) of the system variables $(\Tilde{q},\Tilde{v},\Tilde{u})$. An $m$-tuple of time derivatives $y_{[\Tilde{\kappa}]}=\varphi_{[\Tilde{\kappa}]}(\Tilde{q},\Tilde{v},\Tilde{u},\Tilde{u}_{[1]}, \ldots)$ of this flat output with $\Tilde{\kappa} \leq \Tilde{R}$ and $\# \Tilde{\kappa} = 2(p-k)$ can be introduced as a new input $\Tilde{w}$ by a quasi-static feedback law (\ref{quasistat_fb}) of the generalized state $(\Tilde{q},\Tilde{v})$ if and only if the $2(p-k) \times 2(p-k)$ Jacobian matrix
    \begin{equation}
        \label{jacobian_0}
        \partial_{y_{[0,\Tilde{\kappa}-1]}} 
        \begin{pmatrix}
            F_{\Tilde{q}} \\ F_{\Tilde{v}}
        \end{pmatrix} 
    \end{equation}
    is regular.\footnote{The proof of Proposition \ref{prop_4} is omitted here for conciseness. Interested readers are referred to \cite{HartlExactLinearizationMinimally2023a} for detailed proof of a similar result.}
\end{prop}
Next, we examine all feasible multi-indices $\Tilde{\kappa}$ corresponding to feedback laws that enable rest-to-rest transitions.
\begin{thm}\label{thm_5}
Consider a generalized LC system (\ref{lag_gen_sys_1}) with $m=p-1$ control inputs and a configuration flat output (\ref{flat_output_gen_sys}) as well as the corresponding parameterization of the system variables (\ref{para_map_gen}). Multi-indices $\Tilde{\kappa}$ for which a linear input-output behavior (\ref{input_output_gen}) can be achieved by a quasi-static feedback law (\ref{quasistat_fb}) that is well-defined around equilibrium points can only take one of the two following forms:

i)  For a multi-index $\Tilde{\kappa}=(\Tilde{\kappa}^1, \ldots, \Tilde{\kappa}^{p-1})$ of the form
\begin{equation}
    \label{kappa_1}
    \Tilde{\kappa}^j = 
    \begin{cases}
        2 \text{ for every } j \in \{ 1,\ldots,p-k \}, \\
        0 \text{ for all other } j, %\in \{ p-k+1, \ldots, p-1 \},
    \end{cases}
\end{equation}
there exists a linearizing quasi-static feedback law (\ref{quasistat_fb}) that is well-defined around equilibrium points if and only if there exists a permutation of the components of the flat output 
\begin{equation}\label{thm_5_output_i}
    y = (\underbrace{y^1, \ldots, y^{p-k}}_{\bar{y}}, y^{p-k+1}, \ldots, y^{p-1})
\end{equation}
such that
\begin{equation}\label{thm_5_case_1}
    \Rank{\partial_{\bar{y}}F_{\Tilde{q}}\big|_{y_{s}}}=p-k.
\end{equation}

ii) For a multi-index $\Tilde{\kappa}=(\Tilde{\kappa}^1, \ldots, \Tilde{\kappa}^{p-1})$ of the form
\begin{equation}
    \label{kappa_2}
    \Tilde{\kappa}^j = 
    \begin{cases}
        4 \text{ for one } j \in \{ 1,\ldots,p-k-1 \}, \\
        2 \text{ for all other } j \in \{ 1,\ldots,p-k-1 \}, \\
        0 \text{ for every } j \in \{ p-k, \ldots, p-1 \},
    \end{cases}
\end{equation}
there exists a linearizing quasi-static feedback law (\ref{quasistat_fb}) that is well-defined around equilibrium points if and only if there exists a permutation of the flat output components 
\begin{equation}\label{thm_5_output_ii}
    y = (\underbrace{y^1, \ldots, y^{p-k-1}}_{\bar{y}}, y^{p-k}, \ldots, y^{p-1})
\end{equation}
under which we can find a column $\partial_{\bar{y}^j_{[2]}}\hspace{-1.5mm}F_{\Tilde{q}}|_{y_{s}}$ such that 
\begin{equation}\label{thm_5_case_2}
    \Rank{\partial_{\bar{y}}F_{\Tilde{q}}\big|_{y_{s}}, \hspace{1mm} \partial_{\bar{y}^j_{[2]}}\hspace{-1.5mm}F_{\Tilde{q}}\big|_{y_{s}}} = p-k.
\end{equation}
\end{thm}
\begin{pf}
    According to Proposition \ref{prop_4}, there exists a linearizing feedback law (\ref{quasistat_fb}) that is well-defined around equilibrium points if and only if there exists a submatrix (\ref{jacobian_0}) of the Jacobian
    \begin{equation}
        \label{jacobian_1}
        \partial_{y_{[0, 3]}} 
        \begin{pmatrix}
            F_{\Tilde{q}} \\ F_{\Tilde{v}}
        \end{pmatrix} 
    \end{equation}
    that is regular around $y_s$.  
    Given that (\ref{jacobian_1}) is obtained by omitting $k$ rows from each Jacobian $\partial_{y_{[0, 3]}} F_q$ and $\partial_{y_{[0, 3]}} F_v$ of the parameterization of the classical state, Lemma \ref{lem_2} implies that (\ref{jacobian_1}) evaluated at $y_s$ is of the form
    \begin{equation}
        \label{jacobian_1_ys}
        \left.
        \partial_{y_{[0, 3]}} 
        \begin{pmatrix}
            F_{\Tilde{q}} \\ F_{\Tilde{v}}
        \end{pmatrix} \right|_{y_s} \hspace{-7pt} = \hspace{-3pt}
        \begin{aligned}
        \setlength{\arraycolsep}{0pt}
        \left( \begin{array}{cccc}
            \partial_{y}F_{\Tilde{q}}|_{y_{s}} & 0 & \partial_{y_{[2]}}F_{\Tilde{q}}|_{y_{s}} & 0 \\
            0 & \partial_{y}F_{\Tilde{q}}|_{y_{s}} & 0 & \partial_{y_{[2]}}F_{\Tilde{q}}|_{y_{s}} \\
        \end{array} \right).
        \end{aligned}
    \end{equation}
    Since the $(p-k)\times(p-1)$ matrix $\partial_y F_{\Tilde{q}}$ is obtained by omitting $k$ rows of $\partial_y F_q$ and because of (\ref{lem_1_rank}), we conclude that the rank of $\partial_y F_{\Tilde{q}}|_{y_s}$ is either $p-k$ or \mbox{$p-k-1$} and that $\Rank{\partial_{y_{[2]}}F_{\Tilde{q}}|_{y_s}} \leq 1$. Next, we demonstrate all possibilities of constructing submatrices (\ref{jacobian_0}) that are regular around $y_s$. Considering the special structure of (\ref{jacobian_1_ys}), we begin by selecting $p-k$ or $p-k-1$ linear independent columns of $\partial_y F_{\Tilde{q}}|_{y_s}$. Given $\Rank{\partial_{y_{[2]}}F_{\Tilde{q}}|_{y_s}} \leq 1$, selecting less than $p-k-1$ columns of $\partial_y F_{\Tilde{q}}|_{y_s}$ leads to singular matrices. Thus, constructing all possible submatrices (\ref{jacobian_0}) that are regular around $y_s$ reduces to the following two cases:
    \begin{enumerate}[label=\roman*)]
        \item Given a $\Tilde{\kappa}$ of the form (\ref{kappa_1}), assume the reordered flat output (\ref{thm_5_output_i}) satisfies (\ref{thm_5_case_1}). Since (\ref{conv_1}) implies $y_{[0,\Tilde{\kappa}-1]}=\bar{y}_{[0,1]}$, by Lemma \ref{lem_2}, we obtain a Jacobian of the form
        \begin{equation}
            \label{jacobian_1_i}
            \left. \partial_{y_{[0,\Tilde{\kappa}-1]}} \begin{pmatrix} F_{\Tilde{q}} \\ F_{\Tilde{v}} \end{pmatrix} \right|_{y_{s}} = 
            \begin{aligned}
            \left( \begin{array}{cccc}
                \partial_{\bar{y}}F_{\Tilde{q}}|_{y_{s}} & 0 \\
                0 & \partial_{\bar{y}}F_{\Tilde{q}}|_{y_{s}}  \\
            \end{array} \right).
            \end{aligned} 
        \end{equation}
        Given (\ref{thm_5_case_1}), the $2(p-k)\times2(p-k)$ Jacobian (\ref{jacobian_1_i}) with $\Tilde{\kappa}\leq\Tilde{R}$ and $\#\Tilde{\kappa}=2(p-k)$ is regular around $y_s$.
        \item Given a $\Tilde{\kappa}$ of the form (\ref{kappa_2}), assume there exists a reordered flat output (\ref{thm_5_output_ii}) and a column $\partial_{\bar{y}^j_{[2]}}\hspace{-1.5mm}F_{\Tilde{q}}|_{y_{s}}$ that  satisfy (\ref{thm_5_case_2}). Given (\ref{conv_1}) as well as Lemma \ref{lem_2}, we obtain a Jacobian of the form
        \begin{equation}
            \label{jacobian_1_ys_ii}
            \begin{aligned}
            \partial_{y_{[0,\Tilde{\kappa}-1]}} \hspace{-3pt}
            \begin{pmatrix}
                F_{\Tilde{q}} \\ F_{\Tilde{v}}
            \end{pmatrix}\hspace{-3pt}\Bigg|_{y_{s}} \hspace{-7pt} = \hspace{-3pt}\setlength{\arraycolsep}{0pt}
            \left( \begin{array}{cccc}
                \partial_{\bar{y}}F_{\Tilde{q}}|_{y_{s}} & 0 & \partial_{\bar{y}^j_{[2]}}\hspace{-1.5mm}F_{\Tilde{q}}|_{y_{s}} & 0 \\
                0 & \partial_{\bar{y}}F_{\Tilde{q}}|_{y_{s}} & 0 & \partial_{\bar{y}^j_{[2]}}\hspace{-1.5mm}F_{\Tilde{q}}|_{y_{s}} \\
            \end{array} \right).
            \end{aligned}
        \end{equation}
        Considering the zeros entries of (\ref{jacobian_1_ys_ii}), (\ref{thm_5_case_2}) implies that the $2(p-k) \times 2(p-k)$ submatrix (\ref{jacobian_1_ys_ii}) with $\Tilde{\kappa}\leq\Tilde{R}$ and $\#\Tilde{\kappa}=2(p-k)$ is regular around $y_s$. 
    \end{enumerate}
    Consequently, Proposition \ref{prop_4} implies for each case the existence of a feedback law (\ref{quasistat_fb}) that is well-defined around equilibrium points. Conversely, assume there exists a $\Tilde{\kappa}$ such that (\ref{jacobian_0}) is regular around $y_s$. Since the only possible regular matrices are of the form (\ref{jacobian_1_i}) or (\ref{jacobian_1_ys_ii}), $\Tilde{\kappa}$ can either take the form (\ref{kappa_1}) with a reordered output (\ref{thm_5_output_i}) satisfying (\ref{thm_5_case_1}), or $\Tilde{\kappa}$ can take the form (\ref{kappa_2}) with a reordered output (\ref{thm_5_output_ii}) and a column $\partial_{\bar{y}^j_{[2]}}F_{\Tilde{q}}|_{y_s}$ satisfying (\ref{thm_5_case_2}). \hfill $\Box$
\end{pf}

To formulate a quasi-static feedback law (\ref{quasistat_fb}), we substitute $y_{[\Tilde{\kappa}, \Tilde{R}-1]} = \Tilde{w}_{[0,\Tilde{R}-\Tilde{\kappa}-1]}$ with a feasible $\Tilde{\kappa}$ into the parameterizations (\ref{para_map_gen_q}) of the state $(\Tilde{q}, \Tilde{v})$ and obtain
\begin{equation}
    % \begin{aligned}
    %     \Tilde{q} & = F_{\Tilde{q}}(y_{[0,\Tilde{\kappa}-1]}, \Tilde{w}_{[0,\Tilde{R}-\Tilde{\kappa}-2]}), \\
    %     \Tilde{v} & = F_{\Tilde{v}}(y_{[0,\Tilde{\kappa}-1]}, \Tilde{w}_{[0,\Tilde{R}-\Tilde{\kappa}-1]}).
    % \end{aligned}
    \begin{pmatrix}
        \Tilde{q} \\ \Tilde{v}
    \end{pmatrix} = 
    \begin{pmatrix}
        F_{\Tilde{q}}(y_{[0,\Tilde{\kappa}-1]}, \Tilde{w}_{[0,\Tilde{R}-\Tilde{\kappa}-2]}) \\
        F_{\Tilde{v}}(y_{[0,\Tilde{\kappa}-1]}, \Tilde{w}_{[0,\Tilde{R}-\Tilde{\kappa}-1]})
    \end{pmatrix}.
\end{equation}
Since the multi-index $\Tilde{\kappa}$ is selected such that (\ref{jacobian_0}) is regular around equilibrium points, the implicit function theorem guarantees that a solution of the form
\begin{equation}
    \label{y_kapp}
    y_{[0,\Tilde{\kappa}-1]} = \Tilde{\psi} (\Tilde{q}, \Tilde{v}, \Tilde{w}_{[0,\Tilde{R}-\Tilde{\kappa}-1]})
\end{equation}
exists locally. Subsequently, we can substitute (\ref{y_kapp}) and $y_{[\Tilde{\kappa}, \Tilde{S}]} = \Tilde{w}_{[0,\Tilde{S}-\Tilde{\kappa}]}$ into the flat parameterization (\ref{para_map_gen_u}) of the input $\Tilde{u}$, which further results in the feedback
\begin{equation}
    \Tilde{u} = F_{\Tilde{u}}(\Tilde{\psi} (\Tilde{q}, \Tilde{v}, \Tilde{w}_{[0,\Tilde{R}-\Tilde{\kappa}-1]}), \Tilde{w}_{[0,\Tilde{S}-\Tilde{\kappa}]})
\end{equation}
of the form (\ref{quasistat_fb}).

\subsection{ Exact linearization of the planar aerial manipulator  }

In the following, we illustrate Theorem \ref{thm_5} through a minimally underactuated generalized LC system, namely the planar aerial manipulator presented in Section \ref{sec_lag_con}. Consider the generalized state representation (\ref{plan_manip_gen}) with $(\Tilde{q}, \Tilde{v}) = (x_e, z_e, \theta, v_{x_e}, v_{z_e}, \omega_\theta)$ and $\Tilde{u} = (\phi, F_1, F_2)$. It can be shown that the first three rows of the Jacobian (\ref{jacobian_1_ys}) are given by
\begin{equation}\label{jacobian_case_4}
    \partial_{y_{[0,3]}}F_{\Tilde{q}}\big|_{y_s} = 
    \left(
    \begin{array}{ccc:ccc:ccc:ccc}
        1 & 0 & 0 & 0 & 0 & 0 & \alpha(y^3_s) & 0 & \beta(y^3_s) & 0 & 0 & 0 \\
        0 & 1 & 0 & 0 & 0 & 0 & 0 & 0 & 0 & 0 & 0 & 0 \\
        0 & 0 & 1 & 0 & 0 & 0 & 0 & 0 & 0 & 0 & 0 & 0 \\
    \end{array}
    \right).
\end{equation}
Since $\partial_{y}F_{\Tilde{q}}|_{y_{s}} = I_{3\times3}$ holds, case i) of Theorem \ref{thm_5} applies. Given $p-k=3$, the permutation (\ref{thm_5_output_i}) is clearly given by $y=\bar{y}=(y^1, y^2, y^3)$. Hence, there exists a linearizing quasi-static feedback law (\ref{quasistat_fb}) that is well-defined around equilibrium points. The closed-loop system (\ref{closed_loop}) exhibits a linear input-output behavior (\ref{input_output_gen}) with \mbox{$\Tilde{\kappa}_1=(2,2,2)$}. 

Alternatively, according to case ii), we can select $\bar{y}=(y^2,y^3)$ and have \mbox{$\text{rank}(\partial_{\bar{y}}F_{\Tilde{q}}|_{y_s})=p-k-1=2$}. Next, we choose the column $\partial_{y^3_{[2]}}F_{\Tilde{q}}|_{y_s}$ of (\ref{jacobian_case_4}) such that
\begin{equation*}
    \Rank{\partial_{\bar{y}}F_{\Tilde{q}}\big|_{y_{s}}, \hspace{0mm} \partial_{y^3_{[2]}}\hspace{-1.5mm}F_{\Tilde{q}}\big|_{y_{s}}} = 
    \text{rank}
    \left(
    \begin{array}{cc:c}
        0 & 0 & \beta(y^3_s) \\
        1 & 0 & 0 \\
        0 & 1 & 0 \\
    \end{array}
    \right) \hspace{-0.5mm} = p-k = 3.
\end{equation*}
Consequently, there exists a feedback law (\ref{quasistat_fb}) that is regular around equilibrium points and achieves a linear input-output behavior (\ref{input_output_gen}) with \mbox{$\Tilde{\kappa}_2=(0,2,4)$}.\footnote{ Note that $\Tilde{\kappa}^1_2=0$ implies $y^1=\Tilde{w}^1$, indicating no integrator between the new input $\Tilde{w}^1$ and the flat output $y^1$. }

Let us derive the feedback for $\Tilde{\kappa}_2$. Substituting $y_{[\Tilde{\kappa}_2,\Tilde{R}-1]} = (y^1_{[0,3]}, y^2_{[2,3]}) = ( \Tilde{w}^1_{[0,3]}, \Tilde{w}^2_{[0,1]} )$ into (\ref{plan_manip_flat_para_a})-(\ref{plan_manip_flat_para_b}) yields 
\begin{equation*}\label{eq_sys_fb_2}
    \begin{aligned}
        \Tilde{q} & = F_{\Tilde{q}}( \Tilde{w}^1, \Tilde{w}^1_{[2]}, y^2, \Tilde{w}^2, y^3_{[0,2]} ), \\
        \Tilde{v} & = F_{\Tilde{v}}( \Tilde{w}^1_{[1,3]}, y^2_{[1]}, \Tilde{w}^2_{[0,1]}, y^3_{[0,3]} ),
    \end{aligned}
\end{equation*}
which can locally be solved for $y_{[0, \Tilde{\kappa}_2-1]}=(y^2_{[0,1]}, y^3_{[0,3]})$.\footnote{ Given (\ref{conv_1}), the components of the flat output $y^1$ and $y^3$ do not appear in $y_{[0,\Tilde{\kappa}_2-1]}$ and $y_{[\Tilde{\kappa}_2,\Tilde{R}-1]}$, respectively. } Substituting the solution $y_{[0, \Tilde{\kappa}_2-1]}=\Tilde{\psi}_2(\Tilde{q}, \Tilde{v}, \Tilde{w}^1_{[0,3]}, \Tilde{w}^2_{[0,1]})$ as well as $y_{[\Tilde{\kappa}_2,\Tilde{S}]}=(y^1_{[0,4]}, y^2_{[2,4]}, y^3_{[4]})=(\Tilde{w}^1_{[0,4]}, \Tilde{w}^2_{[0,2]}, \Tilde{w}^3)$ into the parameterization (\ref{plan_manip_flat_para_c}) of the input yields
\begin{equation}\label{quasistat_fb_manip_2}
    \Tilde{u} = \Tilde{\alpha}(x_e, \theta, v_{x_e}, \omega_\theta, \Tilde{w}^1_{[0,4]}, \Tilde{w}^2_{[0,2]}, \Tilde{w}^3),
\end{equation}
representing a feedback law of the desired form (\ref{quasistat_fb}). 

\section{Conclusion}

Our approach illustrates a straightforward method for determining the integrator chain lengths $\Tilde{\kappa}$ of an exactly linearized minimally underactuated generalized LC system with a configuration flat output. The obtained quasi-static feedback (\ref{quasistat_fb}) of the generalized state $(\Tilde{q},\Tilde{v})$ then allows for rest-to-rest transitions and can be systematically derived without requiring the use of successive coordinate transformations or Brunovský states. The presented approach is suitable for any minimally underactuated configuration flat generalized LC system. However, note that neither the existence of a quasi-static feedback law (\ref{quasistat_fb}) that is regular around equilibrium points nor the existence of a quasi-static feedback law (\ref{quasistat_fb}), in general, is guaranteed.
% \fontsize{10}{10.2}\selectfont
\bibliography{references}   
\end{document}